\documentclass[11pt,letterpaper]{article}

\usepackage{amssymb,latexsym,amsmath}
\usepackage[dvipsnames]{color}
\usepackage{graphics}

\definecolor{cite}{named}{Mahogany}
\definecolor{link}{named}{OliveGreen}
\definecolor{anchor}{named}{Gray}

\usepackage[dvips, 
pdftitle={Constructing $\alpha$-critical graphs},
pdfauthor={Carlos E. Valencia},
pdfsubject={Graph Theory},
pdfkeywords={$1$-join, $\alpha$-critical},
colorlinks, 
linkcolor=link, 
citecolor=cite, 
anchorcolor=anchor,
backref, 
pagebackref,
pdfstartview=FitH,
hyperfigures,
bookmarksopenlevel=max, 
]
{hyperref}

\def\qed{\hfill$\Box$}
\def\demo{\noindent{\bf Proof. }}

\def\b1{{\bf 1}}

\newtheorem{theorem}{theorem}[section]

\newtheorem{corollary}[theorem]{Corollary}

\newtheorem{proposition}[theorem]{Proposition}
\newtheorem{claim}[theorem]{Claim}

\begin{document}
\topmargin3mm

\medskip


\vspace{1.5cm}
\begin{center}
{\Large\bf $1$-join composition for $\alpha$-critical graphs.
}
\\
\vspace{3mm}
Carlos E. Valencia\footnote{The firts author was partially supported by SNI.} 
and Marcos I. Leyva\footnote{The second author was partially supported by CONACyT.}    
\\ {\small Departamento de Matem\'aticas} 
\\ {\small Centro de Investigaci\'on y de Estudios Avanzados del IPN} 
\\ {\small Apartado Postal 14--740} 
\\ {\small 07000 M\'exico City, D.F.}
\\ {\small e-mail {\tt cvalenci@math.cinvestav.mx}}
\\
\end{center}

\begin{center}
Version 1, 27 de Julio del 2007.
\end{center}

\medskip


\begin{abstract}
\noindent
Given two graphs $G$ and $H$ its $1$-{\it join} is the graph obtained by taking the disjoint union of $G$ and $H$ and adding all the edges between a nonempty 
subset of vertices of $G$ and a nonempty subset of vertices of $H$.
In general, composition operations of graphs has played a fundamental role in some structural results of graph theory and in particular the $1$-join composition has played an important role in decomposition theorems of several class of graphs such as the claw-free graphs, the bull-free graphs, the perfect graphs, etc. 

A graph $G$ is called {\it $\alpha$-critical} if $\alpha(G\setminus e)> \alpha(G)$ for all the edges $e$ of $G$, where $\alpha(G)$, the {\it stability number} of $G$, 
is equal to the maximum cardinality of a stable set of $G$, and a set of vertices $M$ of $G$ is {\it stable} if no two vertices in $M$ are adjacent.
The study $\alpha$-critical graphs is important, for instance a complete description of $\alpha$-critical graphs would yield a good characterization of the stability number of $G$.

In this paper we give necessary and sufficient conditions that $G$ and $H$ must satisfy in order to its $1$-join will be an $\alpha$-critical graph.
Therefore we get a very useful way to construct basic
$\alpha$-critical graphs using the $1$-join of graphs. 
\end{abstract}

%


\section{Introduction}

One of the most important problems in graph theory consists on obtaining a good characterization of the members of a given class of graphs ${\cal G}$.
Thus, given a family of graphs ${\cal G}$ we would like to have a (de)compostion theorem for ${\cal G}$, that is,
we would like to have a way to construct all the graphs in ${\cal G}$ using some ``basic" class of graphs and some ``basic" operations. 
For instance, if we have a good (de)compostion theorem for ${\cal G}$ we will get an easy way to recognize when a given graph $G$ belongs to ${\cal G}$.
In general, a decomposition theorem has as ingredients some basic subfamilies of ${\cal G}$ and several types of construction of graphs, such as proper $2$-join, balanced skew partitions for Berge graphs and $W$-join, hex-join and generalize $2$-join for claw-free graphs, etc.
Some examples of this decompositions theorems can be seen in~\cite{excluding,claw-free,perfect,cornu1,cornu2}. 

The $1$-join is an important composition operation of graphs that was used to construct decomposition theorems for several important families of graphs.
In general, given two graphs $G$ and $H$ the $1$-join of $G$ and $H$ is the graph obtained by taking the disjoint union of $G$ and $H$ and adding all the edges between a nonempty subset the vertices of $G$ and a nonempty subset of the vertices of $H$.

A subset of vertices $M$ is called a {\it stable} set if any couple of vertices in $M$ are non adjacent.  
The {\it stability number} of a graph $G$ is given by
$$ 
\alpha(G)={\rm max}\{|M| \, | \, M\subset V(G) \mbox{ is a stable set in } G\}. 
$$ 
A graph $G$ is called $\alpha$-critical if $\alpha(G\setminus e)> \alpha(G)$ for all $e\in E(G)$. 
The definition of an $\alpha$-critical graph was firstly introduced in 1949 by Zykov~\cite{Zykov}. 

The $\alpha$-critical graphs have quite interesting properties and have been obtained several beautiful theorems about their structure, see \cite[chapter 12]{LovaszPlummer} for a survey.
The classification and construction of $\alpha$-critical graphs is very important. 
For instance, a complete description of $\alpha$-critical graphs would yield a good characterization of $\alpha(G)$.
The first effort to construct $\alpha$-critical graphs was done by Plummer in \cite{Plummer}; in this article Plummer obtained a family of $\alpha$-critical graphs with an infinite number of elements.

Due to the $\cal NP$-completeness of the problem of calculating $\alpha(G)$ we cannot assume that $\alpha$-critical graphs have a really simple structure,
but several interesting and deep properties has been verified and a certain classification theorems has been proved. 

The {\it defect} $\delta(G)=\vert V(G)\vert-2\alpha(G)=\tau(G)-\alpha(G)$ of a graph plays a central role in the study of $\alpha$-critical graphs.
It was shown in \cite{ErdosGallai} that this defect is non-negative, and the only connected $\alpha$-critical graph with defect zero is ${\cal K}_2$.
One of the most basic facts about $\alpha$-critical graphs is the following theorem done by Hajnal in \cite{Hajnal}:

\begin{theorem}[Hajnal]
If $G$ is a $\alpha$-critical graph, then $\deg(v)\leq \delta(G)+1$ for all $v\in V(G)$.
\end{theorem} 
In particular this theorem implies that the only connected $\alpha$-critical graphs with defect one are the odd cycles.
An odd subdivision of a graph consists in replacing its edges by an odd path.
In this way we can say that the only connected $\alpha$-critical graphs with defect one are the odd subdivision of ${\cal K}_3$.

Andr\'asfai proved in \cite{Andrasfai} a similar result for connected $\alpha$-critical graphs with defect two.
\begin{theorem}[Andr\'asfai]
If $G$ is a $\alpha$-critical graph with defect two, then $G$ is the odd subdivision of ${\cal K}_4$.
\end{theorem}

Sur\'anyi in \cite{Suranyi} gives a characterization of the $\alpha$-critical graph with defect three in the following theorem:
\begin{theorem}[Sur\'anyi]
If $G$ is an $\alpha$-critical graph with defect three, then either $G$ is the odd subdivision of ${\cal K}_5$ or one of the following graphs:
\vspace{22mm} 
\begin{figure}[h] 
\begin{center} 
$ 
\setlength{\unitlength}{.36mm} 
\thicklines\begin{picture}(340,0) 

\hspace{3mm}

\scalebox{0.50}{\includegraphics{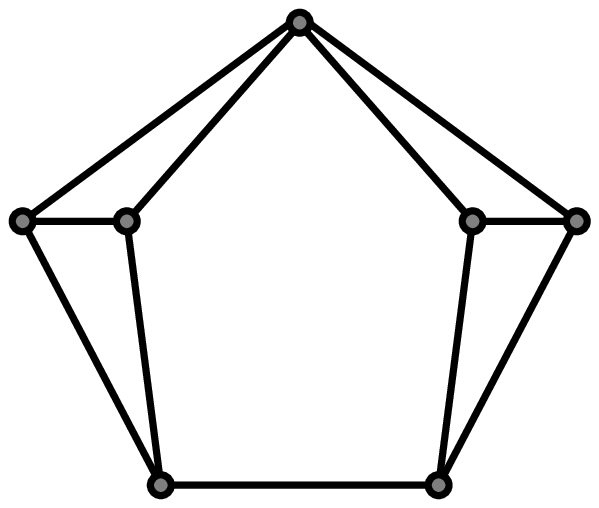}} 

\hspace{14mm}

\scalebox{0.6}{\includegraphics{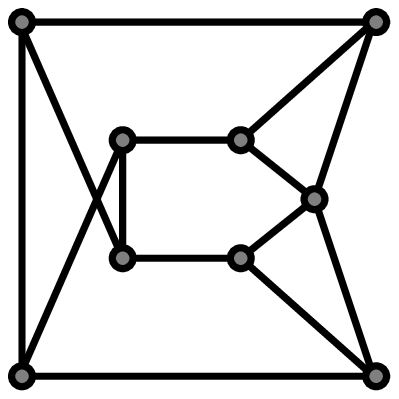}} 

\hspace{14mm}

\scalebox{0.5}{\includegraphics{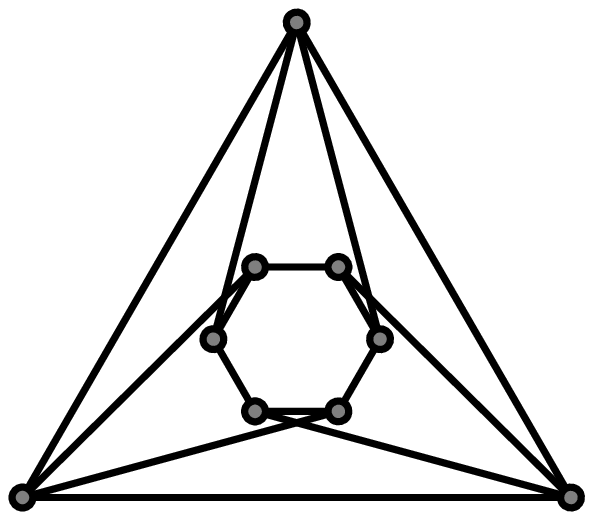}} 

\end{picture} 
$
\vspace{-9mm} 
\end{center} 
\end{figure} 
\end{theorem}

Finally, Lov\'asz in~\cite{Lovasz78} proved that $\alpha$-critical graphs with a fixed defect can be obtained from a finite number of ``basic" graphs by odd subdivision. 

This paper is organized as follows:
In section $2$ we introduce some basic operations of graphs that are useful to construct $\alpha$-critical graphs.
Section $3$ contains the theorem~\ref{1join} that is the main result of the paper.
Theorem~\ref{1join} give necessary and sufficient conditions that $G$ and $H$ must satisfy in order to its $1$-join will be an $\alpha$-critical graph. 
In particular we find that every $\alpha$-critical graph $G$ is the $1$-join of $G\setminus v$ and ${\cal K}_1$.
Moreover, this observation can be use in order to construct $\alpha$-critical graphs in a non trivial way; see~\cite{Delta}.

In order to use the $1$-join composition of graphs we need to have two pairs of graphs $(G, G_0)$ and $(H, H_0)$ where $G_0$ and $H_0$ are maximal induced subgraph of $G$ and $H$ with stability number equal to $\alpha(G_0)=\alpha(G)-1$ and $\alpha(H_0)=\alpha(H)-1$. 
Thus, in section $4$ we characterize all the maximal induced subgraphs $G_0$ of a graph $G$ with stability number equal to $\alpha(G_0)=\alpha(G)-1$ when $G$ is the edge-vertex composition and the $1$-join composition of two graphs.
Finally, in section $5$ we characterize the basic $\alpha$-critical graphs that are the $1$-join of two graphs, that is, we give necessary and sufficient conditions in order that the $1$-join composition of two graph would be splitting free, odd edge-subdivision free and duplication free.

\section{Preliminaries}

Now, we will fix some notation that we will need.
The set of neighborhoods of a subset of vertices $V'\subset V(G)$ 
is equal to 
\[
N_G(V')=\{u\, \vert \, u \mbox{ is adjacent to some } v \in V'\}
\]
and the closed neighborhood of $V'$, denoted by $N_G[V']$, is equal to $N_G[V']=V'\cup N(V')$.
We will denote the set of neighborhoods of $V'$ by $N(V')$ if $G$ is understood for the context. 

\medskip

The induced subgraph of $G$ on a set $V'\subset V(G)$, denoted by $G[V']$,
is the subgraph of $G$ with vertex set equal to $V'$ and edge set equal to 
\[
E(G[V'])=\{e=\{u,v\}\, \vert \, e\in E(G) \mbox{ such that } u,v\in V' \}.
\]

If $V_0\subseteq V(G)$, then $G\setminus V_0$ will denote the induced subgraph $G[V(G)\setminus V_0]$ of $G$.

To get an explanation of terms and symbols see \cite{Diestel}.
\vspace{-4mm}

\subsection{Basic operations on $\alpha$-critical graphs}

In this subsection we will explain with some detail some of the most simple 
operations that preserve the $\alpha$-criticality of the graphs and that 
play an important role in the classification and construction of $\alpha$-critical graphs.
  
\vspace{-4mm}
\paragraph{Edge subdivision}

Let $G$ be a graph and $e=\{u,v\}$ and edge of $G$.
The odd subdivision of $e$ in $G$, denoted by $s(G,e)$, is the new graph given by: 

\begin{description}
\item[$\bullet$] $V(s(G,e))=V(G)\sqcup \{u',v'\}$ and 
\item[$\bullet$] $E(s(G,e))=(E(G)\setminus e)\cup \{u,u'\} \cup \{u',v'\}\cup \{v',v\}$.
\end{description}

\vspace{21mm} 
\begin{figure}[h] 
\begin{center} 
$ 
\setlength{\unitlength}{.4mm} 
\thicklines\begin{picture}(250,0) 

\scalebox{.6}{\includegraphics{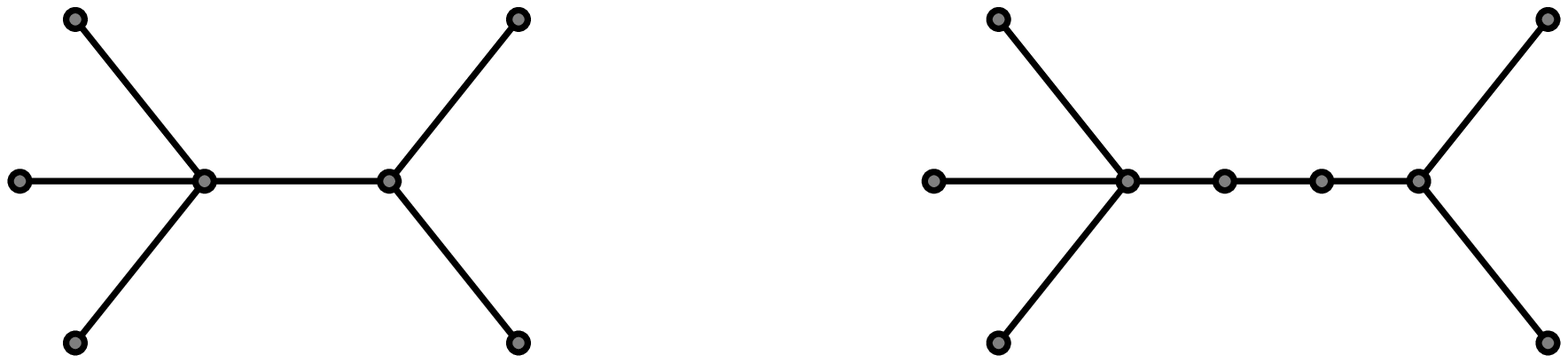}} 

{\LARGE \put(-124,28){$\Longrightarrow$}} 
{\large 
\put(-178,34){$v$}
\put(-205,34){$u$}
\put(-192,33){$e$}
\put(-47,34){$u$}
\put(-33,34){$u'$}
\put(-16,34){$v'$}
\put(-0,34){$v$}
}
\end{picture} 
$ 
\vspace{-18mm}
\end{center} 
\end{figure} 

\newpage

We say that a graph $G$ is {\it odd-subdivision reducible} if there exist a graph $G'$
and an edge $e\in E(G')$ such that $G$ can be obtained by the odd-subdivision of the edge $e$ in $G'$.
A graph $G$ is called {\it odd-subdivision free} if it is not odd-subdivision reducible. 

\vspace{-3mm}

\paragraph{Vertex splitting}

Let $G$ be a graph and $v$ a vertex of $G$. 
The splitting of $v$ in $G$, denoted by $s(G,v)$, is the new graph given by: 
\begin{description}
\item[$\bullet$] $V(s(G,v))=(V(G)\setminus v)\sqcup \{v',v'',u\}$ and 
\item[$\bullet$] $E(s(G,v))=(E(G)\setminus \{\{x,v\}\, \vert \, x\in N(v) \})\cup E_{v'}\cup E_{v''}\cup \{u,v'\} \cup \{u,v''\}$

where $E_{v'}=\{\{v',y\} \, \vert \, y\in N_{v'}\}$, $E_{v''}=\{\{v'',y\} \, \vert \, y\in N_{v''}\}$ 
with $N_{v'}\cup N_{v''}=N_G(v)$, $N_{v'}\cap N_{v''}=\emptyset$ and $N_{v'}, N_{v''}\neq \emptyset$ .
\end{description}

\vspace{21mm} 
\begin{figure}[h] 
\begin{center} 
$ 
\setlength{\unitlength}{.4mm} 
\thicklines\begin{picture}(225,0) 

\scalebox{.6}{\includegraphics{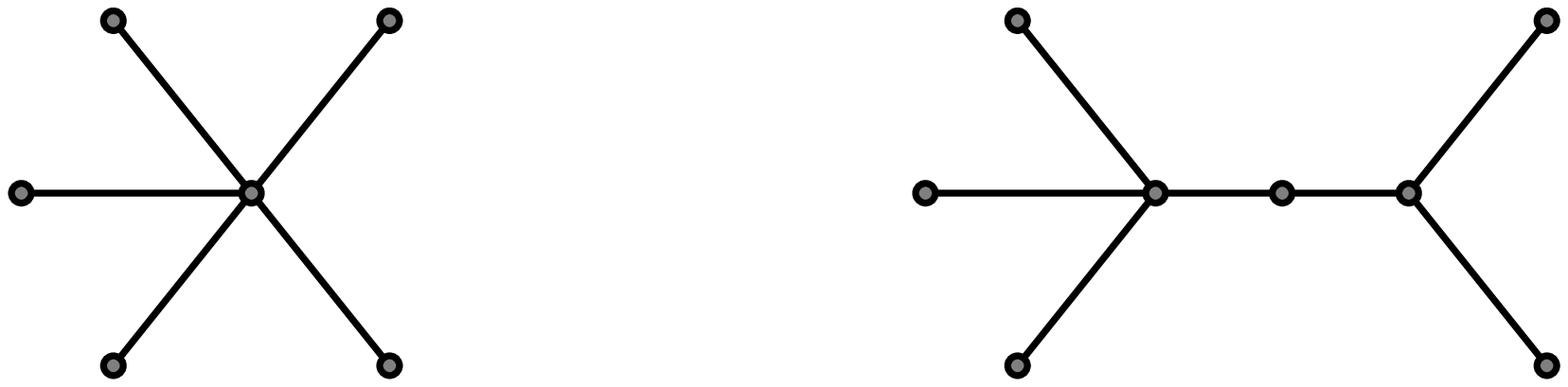}} 
 
{\large 
\put(-115,28){$\Longrightarrow$}
\put(-177,28){$v$}
\put(-37,34){$v'$}
\put(-19,34){$u$}
\put(-1,34){$v''$}
}
\end{picture} 
$ 
\vspace{-5mm}
\end{center} 
\end{figure} 

Note that, if $e=\{v,w\}$ and we take $N_{v'}=N(v)\setminus \{w\}$ and $N_{v''}=\{w\}$, 
then $s(G,v)=s(G,e)$ and therefore the splitting of vertices is a generalization of the odd subdivision of edges.

We say that a graph $G$ is {\it splitting reducible} if there exists a graph $G'$
and a vertex $v\in V(G')$ such that $G$ can be obtained by splitting the vertex $v$ in $G'$.
A graph $G$ is called {\it splitting free} if it is not splitting reducible. 

\vspace{-3mm}

\paragraph{Vertex duplication}
Let $G$ be a graph and $v$ a vertex of $G$, then the duplication of $v$ in $G$, denoted by $d(G,v)$,
is the graph given by:
\begin{description}
\item[$\bullet$] $V(d(G,v))=V(G)\sqcup \{v'\}$ and 
\item[$\bullet$] $E(d(G,v))=E(G)\cup \{\{v',u\} \, \vert \, u\in N[v]\}$.
\end{description}

\vspace{22mm} 
\begin{figure}[h] 
\begin{center} 
$ 
\setlength{\unitlength}{.4mm} 
\thicklines\begin{picture}(190,0) 

\scalebox{.6}{\includegraphics{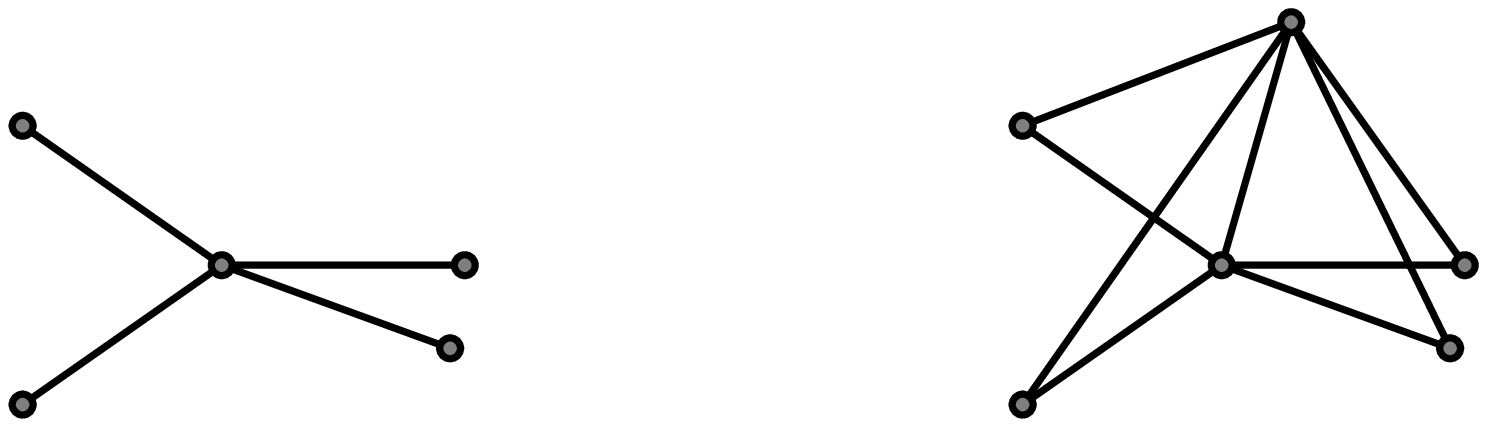}} 
 
{\large 
\put(-110,28){$\Longrightarrow$}
\put(-183,28){$v$}
\put(-16,60){$v'$}
\put(-32,15){$v$}
}
\end{picture} 
$ 
\vspace{-5mm}
\end{center} 
\end{figure} 
 
We say that a graph $G$ is {\it duplication reducible} if there exists a graph $G'$
and a vertex $v\in V(G')$ such that $G$ can be obtained by the duplication of the vertex $v$ in $G'$.
A graph $G$ is called {\it duplication free} if it is not duplication reducible. 

\medskip

Note that ${\cal K}_n=d({\cal K}_{n-1},v)$ for all $v\in V({\cal K}_{n-1})$ 
and $n\in {\mathbb N}$, in particular ${\cal K}_2$ is the duplication 
of the trivial graph ${\cal K}_1$ with only one vertex.

\medskip

Now, we will turn our sight to more general constructions that preserve 
$\alpha$-criticality and permit to construct a huge number of $\alpha$-critical graphs.

\vspace{-3mm}


\paragraph{Edge-Vertex composition} 
Let $G$ and $H$ be $\alpha$-critical graphs, $e=\{v_1,v_2\}$ an edge of $G$ and $v$ a vertex of $H$.
The edge-vertex composition of $(G,e)$ by $(H,v)$, denoted by $c(G,e,H,v)$, is the graph given by
\begin{description}
\item[$\bullet$] $V(c(G,e,H,v))=V(G)\cup (V(H)\setminus v)$
\item[$\bullet$] $E(c(G,e,H,v))=(E(G)\setminus e)\bigcup E(H\setminus v)\bigcup(\{\{v_1,u\}\, \vert \, u\in U_1\}\cup \{\{v_2,u\}\, \vert \, u\in U_2\})$,
 
where $U_1$ and $U_2$ are disjoint non empty sets such that $U_1\cup U_2=N(v)$.
\end{description}

\mbox{}
\vspace{21mm} 
\begin{figure}[h] 
\begin{center} 
$ 
\setlength{\unitlength}{.4mm} 
\thicklines\begin{picture}(257,0) 

\scalebox{.6}{\includegraphics{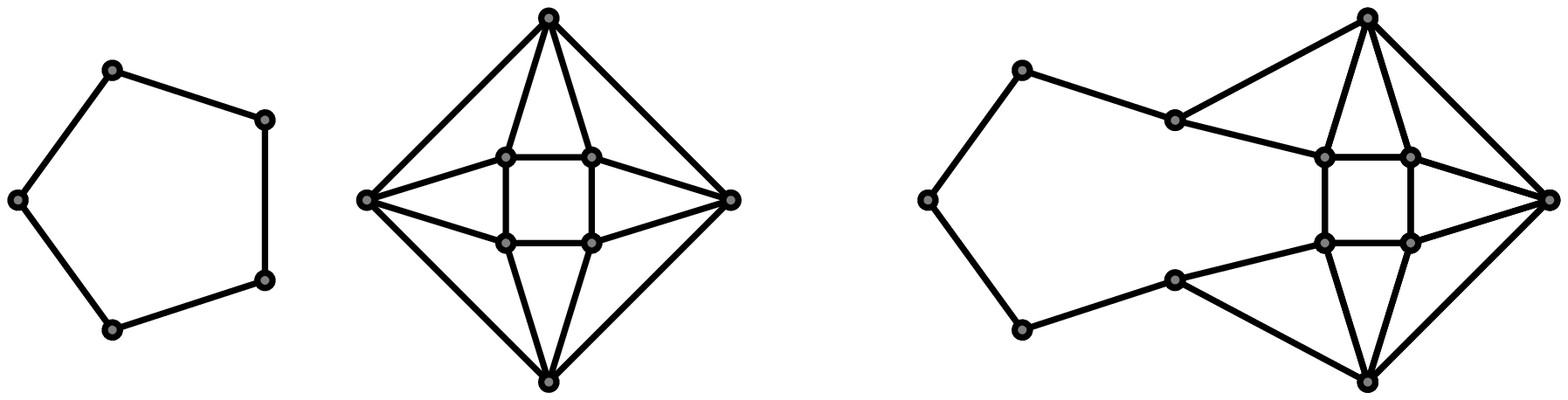}} 
 
{\large 
\put(-133,31){$\Longrightarrow$}
\put(-220,31){$e$}
\put(-212,31){$v$}
\put(-245,30){$G$}
\put(-87,30){$G$}
\put(-178,30){$H$}
\put(-35,30){$H$}
\put(-225,52){$v_1$}
\put(-224,12){$v_2$}
\put(-69,53){$v_1$}
\put(-66,12){$v_2$}
}
\end{picture} 
$ 
\vspace{-6mm}
\end{center} 
\end{figure} 

We say that a graph $W$ is {\it edge-vertex composition reducible} if there exist graphs $G$ and $H$,
$e=\{v_1,v_2\}$ an edge of $G$ and $v$ a vertex of $H$  
such that $W$ is the edge-vertex composition of $(G,e)$ by $(H,v)$.
A graph $W$ is called {\it edge-vertex composition free} if it is not Edge-Vertex composition reducible. 

It is not difficult to see that, $c(G,e,{\cal K}_3,v)=s(G,e)$ and $c({\cal K}_3,e,G,v)=s(G,v)$,
therefore the edge-vertex composition of graphs is a generalization of the odd subdivision of $e$ 
in $G$ and the splitting of $v$ in $G$.

\medskip

The edge-vertex composition introduced by Wessel in~\cite{Wessel} is useful to give a characterization of all 
the $\alpha$-critical graphs with connectivity equal to two, see also Lov\'asz~\cite[Lemma 12.1.5]{Lovasz93}.

\begin{proposition}[\cite{Wessel}]
Let $G$ and $H$ be two $2$-connected $\alpha$-critical graphs, $e$ an edge of $G$ and $v$ a vertex of $H$,
then $c(G,e,H,v)$ is an $\alpha$-critical graph.
Moreover, if $W$ is an $\alpha$-critical graph with connectivity two, then
$W=c(G,e,H,v)$ for some $\alpha$-critical graphs $G$ and $H$.   
\end{proposition}


\section{$1$-join composition for $\alpha$-critical graphs} 

Inspired by a dual form of the edge-vertex composition we obtain a special 
form of the $1$-join composition that is very useful to construct $\alpha$-critical graphs. 
The $1$-join composition has been a very general and useful way to construct graphs 
of several important families of graphs. 
The special form of the $1$-join composition that we will introduce is a very good way of 
constructing $\alpha$-critical graphs. 

\medskip

Let $G$ and $H$ be graphs, $G_0$ and $H_0$ be induced subgraphs of $G$ and $H$ respectively.
The $1$-join composition of $G$ and $H$, denoted by $j(G,G_0,H,H_0)$, is the graph given by:

\begin{description}
\item[$\bullet$] $V((G,G_0,H,H_0))=V(G)\cup V(H)$, 

\vspace{2mm}

\item[$\bullet$] $E(j(G,G_0,H,H_0))=E(G)\bigcup E(H)\bigcup \{\{u,v\}\, \vert \, u\notin V(G_0) \mbox{ and } v\notin V(H_0)\}$.
\end{description}

We say that a graph $J$ is {\it $1$-join reducible} or that admits a {\it $1$-join} if there exist graphs $G$ and $H$ and
induced subgraphs $G_0$ and $H_0$ of $G$ and $H$ respectively,  
such that $J$ is the $1$-join composition of $G$ and $H$.
A graph $J$ is called {\it $1$-join free} if it is not $1$-join reducible.

\vspace{27mm} 
\begin{figure}[h] 
\begin{center} 
$ 
\setlength{\unitlength}{.4mm} 
\thicklines\begin{picture}(273,0) 

\scalebox{.6}{\includegraphics{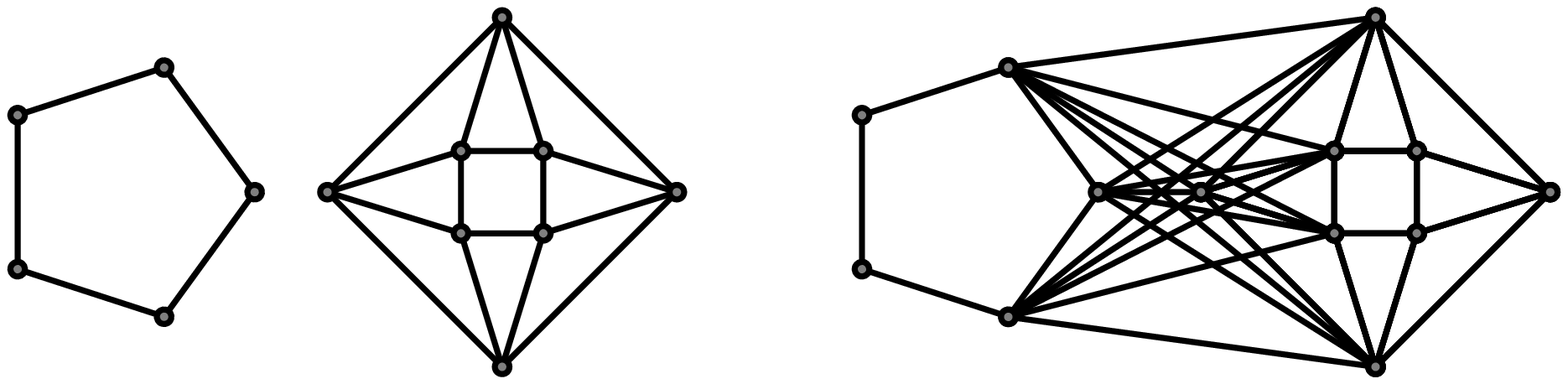}} 
 
{\large 
\put(-155,31){$\Longrightarrow$}
\put(-266,30){$G$}
\put(-297,30){$G_0$}
\put(-199,30){$H$}
\put(-114,30){$G$}
\put(-40,30){$H$}
}
{\footnotesize
\put(-185,32){$H_0$}
}

\end{picture} 
$ 
\vspace{-5mm}
\end{center} 
\end{figure} 

The next theorem gives us the necessary and sufficient conditions in order that
the $1$-join of two graphs to be an $\alpha$-critical graph.

\begin{theorem}\label{1join}
Let $G$ and $H$ be graphs, $G_0$ and $H_0$ be induced subgraphs of $G$ and $H$
with $\alpha(G_0)= \alpha(G)-1$ and $\alpha(H_0)= \alpha(H)-1$, then
\[
\alpha(j(G,G_0,H,H_0))=\alpha(G)+\alpha(H)-1.\eqno(*)
\]
Moreover, $J=j(G,G_0,H,H_0)$ is an $\alpha$-critical graph if and only if
\begin{description}
\item[{\it (i)}] \hspace{1mm} $G_0$ and $H_0$ are maximal (with respect to the inclusion of vertices) induced subgraphs of $G$ and $H$ with $\alpha(G_0)=\alpha(G)-1$ and $\alpha(G_0)=\alpha(G)-1$,

\vspace{2mm}

\item[{\it (ii)}] \hspace{1mm} the edges in $E(G)\setminus E(G_0)$ and $E(H)\setminus E(H_0)$ are $\alpha$-critical edges of $G$ and 
$H$ respectively, 

\vspace{2mm}

\item[{\it (iii)}] \hspace{1mm} if $e\in E(G_0)$, then $e$ is either $\alpha$-critical of $G$ or $G_0$ and 
if $e\in E(H_0)$, then $e$ is either $\alpha$-critical of $H$ or $H_0$.
\end{description}
\end{theorem}
\demo
Let $M$ be a stable set of $J=j(G,G_0,H,H_0)$ and take $M_{G}=M\cap V(G)$ and $M_{H}=M\cap V(H)$.

Firstly, we will prove that $\alpha(J)= \alpha(G)+\alpha(H)-1$.
Since $\alpha(G_0)=\alpha(G)-1$, then there exists a stable set $M_{G_0}$ of $G_0$ with $\vert M_{G_0}\vert=\alpha(G)-1$. 
Let $M'_H$ is a maximum stable set of $H$, then $M_{G_0}\cup M'_H$ is a stable set of $J$.
Thus, clearly $\alpha(j(G,G_0,H,H_0)) \geq \vert M_{G_0}\cup M'_H\vert \geq \alpha(G)+\alpha(H)-1$ because $M_{G_0}\cap M'_H=\emptyset$.

Now, we will prove the other inequality.
Clearly $\vert M_G \vert \leq \alpha(G)$, $\vert M_H \vert \leq \alpha(H)$ and moreover $\alpha(j(G,G_0,H,H_0))\leq \alpha(G)+\alpha(H)$.
Furthermore, if either $\vert M_G\vert\leq \alpha(G)-1$ or $\vert M_H\vert=\alpha(H)-1$, then 
we will have that $\alpha(J)\leq \alpha(G)+\alpha(H)-1$.
Hence it only remains to considerer when either $\vert M_G\vert=\alpha(G)$ or $\vert M_H\vert=\alpha(H)$. 
If we assume that $\vert M_G\vert=\alpha(G)$, then we have that $M_G\cap (V(G)\setminus V(G_0))\neq \emptyset$ because $\alpha(G_0)=\alpha(G)-1$.
Therefore $M_H\subset V(H_0)$ because $N_H(u)=V(H)\setminus V(H_0)$ for all $u\in V(G)\setminus V(G_0)$.
Since $\alpha(H_0)=\alpha(H)-1$ we have that $\vert M_H\vert \leq \alpha(H)-1$ and we obtain that $\alpha(J)\leq \alpha(G)+\alpha(H)-1$.
In the same way if we assume that $\vert M_H\vert=\alpha(H)$ we obtain that $\vert M\vert \leq \alpha(G)+\alpha(H)-1$.
Therefore we can conclude that $\alpha(J)=\vert M\vert \leq \alpha(G)+\alpha(H)-1$.

\medskip

Now, we will prove that $J$ is $\alpha$-critical whenever $G$, $G_0$ and $H$, $H_0$ satisfy the properties $(i),(ii)$ and $(iii)$.
Let us take an edge $e$ of $J$. We will prove that $e$ is a $\alpha$-critical edge.
\begin{description}
\item[$\bullet$] If $e \in E(G)\setminus E(G_0)$, let $M'$ and $M_{H_0}$ be maximum stable sets of $G\setminus e$ and $H_0$ respectively.
Then $M=M'\cup M_{H_0}$ is  a stable set of $J\setminus e$ with $(\alpha(G)+1)+(\alpha(H)-1)=\alpha(G)+\alpha(H)$ vertices
and therefore $e$ is an $\alpha$-critical edge of $J$.

\vspace{2mm}

\item[$\bullet$] If $e \in E(G_0)$ and $e$ is an $\alpha$-critical edge of $G$, then we can use the same argument that in the previous case.

\vspace{2mm}

\item[$\bullet$] If $e \in E(G_0)$ and $e$ is an $\alpha$-critical edge of $G_0$, let $M'_{G_0}$ and $M_H$ be maximum stable sets of $G_0\setminus e$ and $H$ respectively.
Then $M=M'_{G_0}\cup M_H$ is  a stable set of $J\setminus e$ with $\alpha(G)+\alpha(H)$ vertices
and therefore $e$ is an $\alpha$-critical edge of $J$.

\vspace{2mm}

\item[$\bullet$] By symmetry of $J$ with respect to $G$ and $H$ we can use the same arguments to prove that the edges of $H$ are $\alpha$-critical edges of $J$.

\vspace{2mm}

\item[$\bullet$] If $e=\{u,v\}$ for some $u\notin V(G_0)$ and $v\notin V(H_0)$.
Since $G_0$ and $H_0$ are maximal induced subgraphs of $G$ and $H$, then
there exist $M_u$ and $M_v$ stable set of $G$ and $H$ respectively with $u\in M_u$, 
$v\in M_v$, $\vert M_u\vert=\alpha(G)$, $\vert M_v\vert=\alpha(H)$, 
$M_u\setminus u \subset V(G_0)$ and $M_v\setminus v \subset V(H_0)$. 
Hence $M=M_u\cup M_v$ is a stable set of $J\setminus e$
with $\alpha(G)+\alpha(H)$ vertices and therefore $e$ is an $\alpha$-critical edge of $J$.

\end{description}

\medskip

To finish we will prove that if $j(G,G_0,H,H_0)$ is an $\alpha$-critical graph, 
then $G$, $G_0$ and $H$, $H_0$ satisfy the properties $(i),(ii)$ and $(iii)$.

Let us assume that neither $G_0$ is a maximal induced subgraph of $G$ nor $H_0$ is a maximal induced subgraph of $H$.
Hence there exist $G'_0$ and $H'_0$ maximal induced subgraphs of $G$ and $H$ with 
$V(G_0)\subseteq V(G'_0)$ and $V(H_0)\subseteq V(H'_0)$.
Using $(*)$ we have that $j(G,G'_0,H,H'_0)$ is a graph with 
stability number equal to $\alpha(G)+\alpha(H)-1$ and since $j(G,G'_0,H,H'_0)$ is a 
spanning subgraph of $J$ with $E(j(G,G'_0,H,H'_0)) \subsetneq E(J)$, 
then $J$ can not be an $\alpha$-critical graph; a contradiction.

Now, let us take an edge $e$ of $G$, if $e\in E(G)\setminus E(G_0)$ and $e$ is not an $\alpha$-critical edge of $G$, 
then $\alpha(G\setminus e)=\alpha(G)$ and clearly $G_0$ is an induced subgraph of $G\setminus e$ with $\alpha(G_0)=\alpha(G\setminus e)$.
Since $j(G\setminus e, G_0,H,H_0)=J\setminus e$ and using $(*)$ we obtain that $\alpha(j(G\setminus e,G_0,H,H_0))=\alpha(J)$;
a contradiction to the fact that $J$ is an $\alpha$-critical graph. 
In a similar way if $e\in E(G_0)$ and $e$ is not an $\alpha$-critical edge of $G$ or $G_0$, 
then $\alpha(G\setminus e)=\alpha(G)$ and $G_0\setminus e$ is an induced subgraph 
of $G\setminus e$ with $\alpha(G_0\setminus e)=\alpha(G\setminus e)-1$.
Hence using $(*)$ we obtain that $\alpha(j(G\setminus e,G_0\setminus e,H,H_0))=\alpha(J)$ 
and therefore $J$ would not be an $\alpha$-critical graph; a contradiction. 

By the symmetry of $j(G,G_0,H,H_0)$ with respect to $G$ and $H$, we can apply the same arguments when $e$ is an edge of $H$.
\qed

\medskip 

If we take $G={\cal K}_1$, $G_0$ as the empty graph (note that $G_0=G\setminus N[u]$ 
with $u$ the only one vertex of $G$) and $H_0=H\setminus N[v]$ for some vertex of $H$, 
then $j(G,G_0,H,H_0)=d(H,v)$.
Therefore we can think the $1$-join of graphs as a generalization of the duplication 
of a vertex.

Note that in contrast to the edge-vertex composition of graphs, the $1$-join of two graphs
allows to construct $\alpha$-critical graphs with high connectivity.

The previous theorem tells us that the $1$-join of two graphs $G$ and $H$ can be an $\alpha$-critical even if
$G$ and $H$ are not $\alpha$-critical graphs.
In fact, if either $G$ is an $\alpha$-critical graph or $G_0$ is an $\alpha$-critical 
graph and the edges in $E(G)\setminus E(G_0)$ are $\alpha$-critical edges of $G$, and 
either $H$ is an $\alpha$-critical graph or $H_0$ is $\alpha$-critical graph and the 
edges in $E(H)\setminus E(H_0)$ are $\alpha$-critical edges of $H$, then $j(G,G_0,H,H_0)$ 
is an $\alpha$-critical graph.

\medskip

An implication of the theorem~\ref{1join} is the following corollary:
\begin{corollary}
Let $G$ and $H$ be graphs, $G_0$ and $H_0$ be induced subgraphs of $G$ and $H$,
then $J=j(G,G_0,H,H_0)$ is and $\alpha$-critical graph if and only if $G_1=j(G,G_0,{\cal K}_1,\emptyset)$ 
and $H_1=j(H,H_0,{\cal K}_1,\emptyset)$ are $\alpha$-critical graphs.
\end{corollary}
\demo
It follows directly from theorem~\ref{1join}.
\qed

See~\cite{bixby} for a similar result for perfect graphs. 

\medskip

Now, we will study some special cases of the $1$-join of graphs.
\begin{corollary}
Let $G$ and $H$ be graphs, $G_0$ and $H_0$ be induced subgraphs of $G$ and $H$
with $\alpha(G_0)= \alpha(G)-1$ and $\alpha(H_0)= \alpha(H)-1$ such that $G$, 
$G_0$, $H$ and $H_0$ satisfy the conditions $(i)$, $(ii)$ and $(iii)$,
then
\begin{description}
\item[$\bullet$] $j(G,G_0,H,H_0)\setminus u$ is an $\alpha$-critical graph if and only if $G_0=G\setminus u$ for some $u\in V(G)$.

\vspace{2mm}

\item[$\bullet$] $j(G,G_0,H,H_0)\setminus \{u,v\}$ is an $\alpha$-critical graph if and only if $G_0=G\setminus u$ for some $u\in V(G)$ and $H_0=H\setminus v$ for some $v\in V(H)$.
\end{description}
\end{corollary}
\demo
The result is followed by using the same arguments as in the 
proof of the previous theorem and the observation that
if $G_0=G\setminus N[v]$ for some $v\in V(G)$, $e\in E(G_0)$ 
is $\alpha$-critical graph of $G$ but not an $\alpha$-critical 
graph of $G_0$ and $M$ is a maximum stable set of $G\setminus e$,
then $M\cap \{v\}=\emptyset$.
\qed

\medskip

Using two graphs $G$ and $H$ as blocks and the theorem $4$ and corollary $2$
we can construct $\alpha$-critical graphs with defect equal to $\delta(G)+\delta(H)+2$,
$\delta(G)+\delta(H)+1$ and $\delta(G)+\delta(H)$.

The following graph is an example of an $\alpha$-critical graph $j(G,G_0,H,H_0)$
with $G={\cal K}_1$, $V(G_0)=\emptyset$ and $H$ is not an $\alpha$-critical graph. 

\vspace{29mm}
\begin{figure}[h] 
\begin{center} 
$ 
\setlength{\unitlength}{.4mm} 
\thicklines\begin{picture}(130,0) 

\scalebox{0.65}{\includegraphics{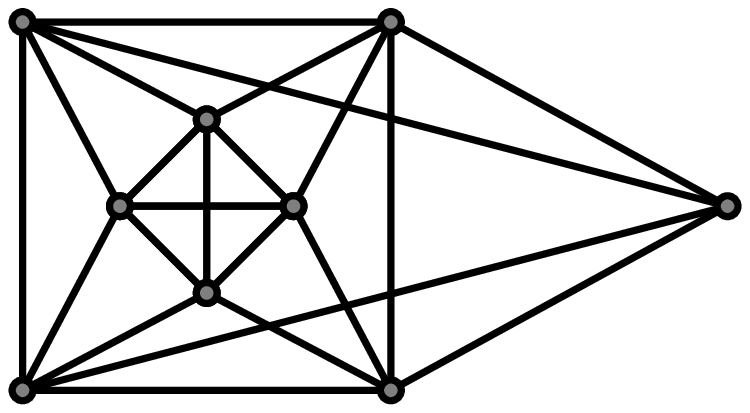}} 
 
{\normalsize 
\put(-19,44){$v$}
\put(-110,73){$H$}
}
{\small
\put(-124,50){$H_0$}
}
\end{picture} 
$ 
\vspace{-10mm}
\end{center} 
\end{figure}
 
\medskip

Theorem~\ref{1join}  tell us that in order to construct $\alpha$-critical 
graphs from two $\alpha$-critical graphs $G$ and $H$, then we need
to have maximal induced subgraphs $G_0$ and $H_0$ of $G$ and $H$ respectively.
The following theorem gives us a good source of this type of induced graphs.

\begin{theorem}\label{max}
Let $G$ be a graph and $v$ a vertex of $G$. 
Then every edge of $G$ incident with $v$ is an $\alpha$-critical edge of $G$ if and only if
$G_v=G\setminus N[v]$ is a maximal induced subgraph of $G$ with $\alpha(G_v)=\alpha(G)-1$. 
\end{theorem}
\demo
$(\Rightarrow)$
Let $M$ be a stable set of $G_v=G\setminus N[v]$. 
We have that $\vert M \vert \leq \alpha(G)-1$.
Since if we suppose that $\vert M \vert > \alpha(G)-1$, then $M\cup \{v\}$ would be a stable set of $G$ with 
$\vert M\cup \{v\} \vert > \alpha(G)$; a contradiction.

Now, let $w\in N(v)$. Since $e=\{v,w\}$ is an $\alpha$-critical edge of $G$, 
then there exists a stable set $M_w$ of $G\setminus e$ with $\vert M_w \vert =\alpha(G)+1$ and $v,w \in M_w$.
Hence $\alpha(G_v)=\alpha(G)-1$ because $M_w\setminus \{v,w\}$ is a stable
set of $G_v$ ($M_w \cap N[v]=\{v,w\}$) with $\alpha(G)-1$ vertices.
Moreover, $G_v$ is maximal because for all $w\in N(v)$ we have that $M_w\setminus v$ is a stable set of 
$G[V(G_v)\cup \{w\}]$ with $\alpha(G)$ vertices and 
$M_w\setminus w$ is a stable set of $G[V(G_v)\cup \{v\}]$ with $\alpha(G)$ vertices.

$(\Leftarrow)$
Let  $w\in N(v)$ and $e=\{v,w\}$, since $G_v$ is a maximal induced subgraph of $G$ with $\alpha(G_v)=\alpha(G)-1$,
then $\alpha(G[(V(G)\setminus N[v])\cup \{w\}])=\alpha(G)$.
Furthermore, if $M$ a maximum stable set of $G[(V(G)\setminus N[v])\cup \{w\}]$,
then $M\cup \{v\}$ is a stable set of $G\setminus e$ and therefore $e$ is an $\alpha$-critical edge of $G$. 
\qed

\begin{corollary}\label{maxall}
A graph $G$ is $\alpha$-critical if and only if $G_v=G\setminus N[v]$ is a maximal 
induced subgraph of $G$ with $\alpha(G_v)=\alpha(G)-1$ for all $v\in V(G)$.
\end{corollary}
\demo
It follows directly from theorem~\ref{max}.
\qed

\medskip

By theorems~\ref{1join} and \ref{max} and since $\alpha(G\setminus v)=\alpha(G)$, then we have that every $\alpha$-critical graph is the $1$-join of the pairs $(G\setminus v,G\setminus N[v])$ and $({\cal K}_1,\emptyset)$.
This trivial way to compose an $\alpha$-critical graph as the $1$-join of two graphs is interesting, 
for instance in~\cite{Delta} this observation is used in order to study the special case when $\alpha(G)=2$.
The case when $\alpha(G)=2$ is the most simple case because we have that $G_v$ is a complete graph.
More precisely, in~\cite{Delta} 
it is described a process that uses a $\Delta$-hypergraph $H$ 
as a base in order to construct a new $\alpha$-critical graph $H'$ with $\alpha(H')=2$. 

The construction described in~\cite{Delta} is a generalization of the dual (under taking the complement of graphs) of Mycielski construction of triangle-free graphs, 
see~\cite{Mycielski} or \cite[pag. 253]{Kocay} for a description of Mycielski construction.

Let $a$ be a fixed natural number, for $a=2$ and $a=3$ the 
corollary~\ref{maxall} can be use
in order to construct a polynomial algorithm (with complexity $\vert V(G)\vert^a$) that recognizes $\alpha$-critical graphs
with stability number equal to $a$, see~\cite{Marcos} for details of the algorithm for graphs with stability number $2$ and $3$.

\medskip

Note that not all the maximal induced subgraphs $G_0$ of an $\alpha$-critical graph $G$
with $\alpha(G_0)=\alpha(G)-1$ are equal to $G_v=G\setminus N[v]$ for some $v\in V(G)$.
For instance in the following two graphs take $H_0$ as the induced subgraph by $V'=\{v_1,v_2,v_3,v_4,v_5\}$. 
Is not difficult to check that $H_0$ is a maximal induced subgraph with $\alpha(H_0)=\alpha(H)-1$ but $H_0\neq H_v$ for all $v\in V(H)$.

\vspace{33mm}
\begin{figure}[h] 
\begin{center} 
$ 
\setlength{\unitlength}{.4mm} 
\thicklines\begin{picture}(220,0) 

\scalebox{0.75}{\includegraphics{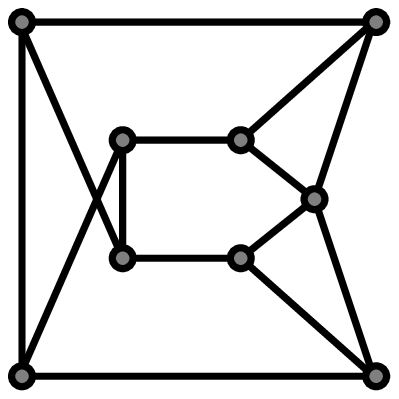}} 
 
{\normalsize 
\put(-12,34){$v_1$}
\put(-36,53){$v_2$}
\put(-55,53){$v_3$}
\put(-55,18){$v_4$}
\put(-35,18){$v_5$}
\put(-41,34){$H_0$}
\put(0,-1){$v_6$}
\put(-83,-1){$v_7$}
\put(-83,70){$v_8$}
\put(0,70){$v_9$}
\put(-45,59){$H$}
}

\hspace{25mm}
\scalebox{0.65}{\includegraphics{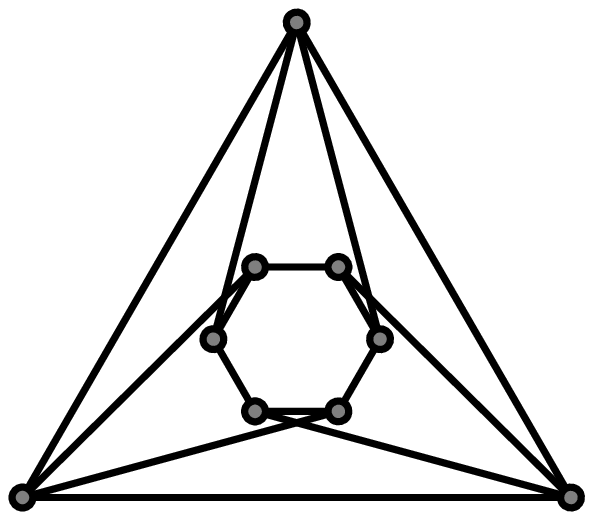}} 
{\normalsize 
\put(-104,-1){$v_1$}
\put(0,-1){$v_2$}
\put(-45,80){$v_3$}
\put(-39,14){$v_5$}
\put(-34,22){$v_4$}
\put(-30,12){$H_0$}
\put(-66,15){$v_6$}
\put(-69,22){$v_7$}
\put(-56,44){$v_8$}
\put(-48,45){$v_9$}
\put(-70,36){$H$}
}

\end{picture} 
$ 
\vspace{-5mm}
\end{center} 
\end{figure} 

In the first graph, $H_0$ is an $\alpha$-critical graph (isomorphic to $C_5$) 
and in the second graph $H_0$ is not an $\alpha$-critical graph (because $\{v_1,v_3\}$ is not an $\alpha$-critical edge of $G$).

A maximal induced subgraph $G_0$ of a graph $G$ with $\alpha(G_0)=\alpha(G)-1$ is called canonical 
if $G_0=G\setminus N[v]$ for some $v\in V(G)$ and is called non-canonical if it is not canonical.


\section{Maximal induced subgraphs} 

In order to use theorem~\ref{1join} to construct $\alpha$-critical graphs we need to have two pairs of graphs $(G, G_0)$ and $(H, H_0)$ where $G_0$ and $H_0$ are maximal induced subgraph of $G$ and $H$ respectively with stability number equal to $\alpha(G_0)=\alpha(G)-1$ and $\alpha(H_0)=\alpha(H)-1$ 
and such that it satisfies $(i), (ii)$ and $(iii)$ in theorem~\ref{1join}.
If $G$ is $\alpha$-critical, then we only need that $G_0$ will be a maximal 
induced subgraph of $G$ with $\alpha(G_0)=\alpha(G)-1$.
The next two theorems give us a characterizations of the maximal induced subgraphs $G_0$ of $G$
with $\alpha(G_0)=\alpha(G)-1$ when $G$ is the edge-vertex composition and the $1$-join of two graphs.

\medskip

First we solve the case when $G$ is the edge-vertex composition of two graphs.

\begin{theorem}\label{maximalcomposition}
Let $G$ and $H$ be $\alpha$-critical graphs, $e=\{v_1,v_2\}$ an edge of $G$ and $v$ a vertex of $H$,
then $W_0$ is a maximal induced subgraph of $W=c(G,e,H,v)$ with $\alpha(W_0)=\alpha(W)-1$ if and only if 
$W_0$ is one of the following graphs:
\begin{description}
\item[(i)] \hspace{1mm} $W[V(G'_0)\sqcup (V(H)\setminus v)]$ where $G'_0$ a maximal induced subgraph of 
$G$ with $\alpha(G'_0)= \alpha(G)-1$.

\vspace{2mm}

\item[(ii)] \hspace{1mm} $c(G,e,H'_0,v)=W[V(G) \sqcup (V(H'_0)\setminus v)]$ where $v\in V(H'_0)$ and $H'_0$ is a maximal 
induced subgraph of $H$ with $\alpha(H'_0)= \alpha(H)-1$.

\vspace{2mm}
 
\item[(iii)] \hspace{1mm} $W[(V(G) \setminus v_i)\sqcup V(H'_0)]$  for $i=1,2$  and where $v\notin V(H'_0)$ and $H'_0$ is a maximal induced 
subgraph of $H$ with $\alpha(H'_0)= \alpha(H)-1$. 

\vspace{2mm}
 
\item[(iv)] \hspace{1mm} $W[V(G'_0)\cup v_i \sqcup (V(H)\setminus (N_W(v_j)\cup v))]$ for $i,j\in \{1,2\}$ and $i\neq j$ 
where $v_i\notin V(G'_0)$ and $G'_0$ a maximal induced subgraph of $G$ with $\alpha(G'_0)= \alpha(G)-1$.  
\end{description}
\end{theorem}
\demo
Before we start the proof we will prove the following basic result:

\medskip

\begin{claim} 
Let $G$ and $H$ be graphs, $e=\{v_1,v_2\}$ an edge of $G$ and $v$ a vertex of $H$, then
\[
\alpha(c(G,e,H,v))\leq \alpha(G)+\alpha(H).
\]
\end{claim}
\demo
Let $M$ be a maximum stable set of $c(G,e,H,v)$ and let us take $M_G=V(G)\cap M$ and $M_H=V(H)\cap M$.
Clearly $\vert M_G\vert \leq \alpha(G)+1$ and $\vert M_H\vert \leq \alpha(H)$ because
we have that $M_G$ and $M_H$ are stable sets of $G\setminus e$ and $H\setminus v$ respectively.

If $\{v_1,v_2\}\not\subset M$, then we have that $\vert M_G\vert \leq \alpha(G)$ because in this case $M_G$ would be a stable set of $G$.
Hence we have that if $\vert M_G \vert =\alpha(G)+1$, then $v_1,v_2\in M$.
Also, if $v_1,v_2\in M$, then we will have that $M_H\cup\{v\}$ is a stable set of $H$ and therefore $\vert M_H\vert \leq \alpha(H)-1$. 
Using this two previous observation is not difficult to conclude that
\[
\alpha(c(G,e,H,v))=\vert M\vert =\vert M_G\vert+\vert M_H\vert\leq \alpha(G)+\alpha(H).\vspace{-8.8mm}
\]
\qed

\vspace{4mm}

\noindent $(\Leftarrow)$
$(i)$ Let $W_0=W[V(G'_0)\sqcup (V(H)\setminus v)]$ where $G'_0$ is a maximal 
induced subgraph of $G$ with $\alpha(G'_0)= \alpha(G)-1$.
Since $\alpha(H\setminus v)=\alpha(H)$ ($H$ is $\alpha$-critical), then using
claim 1 
and the fact that the stability number of the disjoint
union of two graphs is equal to the sum of the stability numbers of this two graphs,
is not difficult to see that $\alpha(W_0)\leq \alpha(G)-1+\alpha(H)=\alpha(W)-1$.

On the other hand, for $i=1,2$ let us take $M_{v_i}=(M\setminus v)\cup \{v_i\}$ where $M$ 
is a maximum stable set of $H\setminus e'$ with $e'=\{v,w\}\in E(H)$ and $w\in N_W(v_j)\cap V(H)$. 
Clearly $M_{v_1}$ and $M_{v_2}$ are stable sets of $H'=W[(V(H)\setminus v)\cup \{v_1,v_2\}]$ with $\alpha(H)+1$ vertices.

Let $M_{G_0}$ be a maximum stable set of $G'_0$, since $v_1\notin M_{v_2}$, $v_2\notin M_{v_1}$ and $\{v_1,v_2\}\not\subset M_{G'_0}$, then we can conclude that
\[
\alpha(W_0)=\alpha(G'_0)+\alpha(H\setminus v)=\alpha(G)-1+\alpha(H)=\alpha(W)-1.
\]
It only remains to prove that $W_0$ is a maximal induced subgraph of $W$ with $\alpha(W_0)=\alpha(W)-1$.
Let us take $u\in V(W)\setminus V(W_0)=V(G)\setminus V(G_0)$.
By the maximality of $G_0$ we have that there exists a maximum stable set $M_u=M_{G'_0}\cup\{u\}$ 
(where $M_{G'_0}$ is a maximum stable set of $G'_0$) of $G''_0=G[V(G'_0)\cup \{u\}]$
with $\{v_1,v_2\}\not\subset M_u$ and $\alpha(G)$ vertices.
Since $\alpha(H\setminus v)=\alpha(H)$, then there exists a stable set $S$ of $H\setminus v$ with $\alpha(H)$ vertices.
If $v_1,v_2\notin M_{u}$, then $M_u\cup S$ is a maximum stable set of $W[V(W_0)\cup \{u\}]$ with $\alpha(W_0)+1=\alpha(G)+\alpha(H)$ vertices.
If $v_i\in M_{u}$ and $v_j\notin M_{u}$ for $\{i,j\}=\{1,2\}$, then $M_u\cup M_{v_i}$ is a maximum stable set of $W[V(W_0)\cup \{u\}]$ with $\alpha(W_0)+1=\alpha(G)+\alpha(H)$ vertices. 
Therefore $W_0$ is maximal.

\medskip

$(ii)$ Let $W_0=c(G,e,H'_0,v)$ where  $v\in V(H'_0)$ and $H'_0$ is a maximal induced subgraph of $H$ with $\alpha(H'_0)= \alpha(H)-1$. 
Let $M_{H'_0}$ be a maximum stable set of $H'_0$, $M=M'\setminus v_1$ with $M'$ the maximum stable set of $G\setminus \{v_1,u\}$ 
where $u\in N(v_1)\setminus v_2$ and $M_{v_1,v_2}$ the maximum stable set of 
$G\setminus e$ with $v_1,v_2\in M_{v_1,v_2}$ and $\vert M_{v_1,v_2} \vert =\alpha(G)+1$.

If $v\in M_{H'_0}$, then $(M_{H'_0}\setminus v)\cup M_{v_1,v_2}$ is a stable set of $c(G,e,H_0,v)$ 
with $\alpha(G)+\alpha(H)-1$ vertices and 
if $v\notin M_{H'_0}$, then $M_{H'_0}\cup M$ is a stable set of $c(G,e,H_0,v)$ ($v_1,v_2\notin M$) with $\alpha(G)+\alpha(H)-1$ 
vertices.

Therefore using claim 1 
we can conclude that $\alpha(c(G,e,H'_0,v))=\alpha(G)+\alpha(H)-1$. 
Moreover using a similar argument is not difficult to prove that if $H'_0$ is a maximal induced 
subgraph of $H$ with $\alpha(H'_0)=\alpha(H)-1$, then $c(G,e,H'_0,v)$ is a maximal
induced subgraph of $c(G,e,H,v)$ with $\alpha(c(G,e,H'_0,v))=\alpha(G)+\alpha(H)-1$.

\medskip

$(iii)$
Let $W_0=W[(V(G) \setminus v_1)\sqcup V(H'_0)]$ with $v\notin V(H'_0)$ and $H'_0$ be a maximal induced 
subgraph of $H$ with $\alpha(H'_0)= \alpha(H)-1$.
Let $M_v$ be a maximum stable set of $H[V(H'_0)\cup\{v\}]$ with $v\in M_v$, $M_v\setminus v\subset V(H'_0)$ 
and $\vert M_v\vert = \alpha(H)$.
Hence $(M_v\setminus v)\cup M'_{v_2}$, where $v_2 \in M'_{v_2}=M_{v_1,v_2}\setminus v_1$ is a maximum stable set of 
$G\setminus v_1$ with $\alpha(G)$ vertices, is a stable set of $W[V(G\setminus v_1)\cup V(H'_0)]$ with $\alpha(G)+\alpha(H)-1$ vertices.
Moreover, since $\alpha(W[V(G\setminus v_1)\cup V(H'_0)])\leq \alpha(G)+\alpha(H'_0)=\alpha(G)+\alpha(H)-1$, then
$\alpha(W[V(G\setminus v_1)\cup V(H'_0)])=\alpha(G)+\alpha(H)-1$.

The induced subgraph $W[V(G\setminus v_1)\cup V(H'_0)]$ is maximal, since if $u\in V(H)\setminus (V(H'_0)\cup v)$,
then $(M_{v_1,v_2}\setminus v_2)\cup M_u$, where $M_u$ is a stable set of $H[V(H'_0)\cup\{u\}]$ with $\alpha(H)$ vertices,
is a stable set of $W[V(G\setminus v_1)\cup V(H'_0)\cup \{u\}]$ with $\alpha(G)+\alpha(H)$ vertices 
and if $u=v_1$, then $M_{v_1,v_2}\cup (M_v\setminus v)$ is a stable set of $W[V(G)\cup V(H'_0)]$ with $\alpha(G)+\alpha(H)$ vertices. 

\medskip

We can use exactly the same arguments when $W_0=W[V(G\setminus v_2)\cup V(H'_0)]$.

\medskip

$(iv)$
Let $W_0=W[V(G'_0)\cup v_1 \sqcup (V(H)\setminus (N_W(v_2)\cup v))]$ with 
$v_1\notin V(G'_0)$ and $G'_0$ a maximal induced subgraph of $G$ with $\alpha(G'_0)= \alpha(G)-1$.
Clearly $\alpha(W_0)\leq \alpha(G'_0)+\alpha(W[\{v_1\}\cup(V(H)\setminus (N_W(v_2)\cup v))])=\alpha(G)-1+\alpha(H)$. 

Since $G'_0$ is a maximal induced subgraph of $G$ with $\alpha(G'_0)=\alpha(G)-1$,
then for all $u \notin V(G'_0)$ there exists a maximum stable set $M_u$ of $G$ with $u\in M_u$.
Let $N_H$ be a maximum stable set of $H$. If $v\notin N_H$, then $N_W(N_H)\cap \{v_1,v_2\}\neq \emptyset$,
more precisely if $N_W(v_j)\cap N_H=\emptyset$, then $v_i\in N_W(N_H)$ for $i\neq j$. 

For $k=1,2$, take $N_{v_k}=N\setminus v$, where $N$ is a maximum stable set of $H\setminus \{v,w\}$ with $w \in N_W(v_k)\cap V(H)$.
Since $(N_W(v_k)\cup v)\cap N_{v_k}=\emptyset$, then we have that $M_{v_1}\cup N_{v_2}$ ($M_{v_i}$ is defined in $(i)$)is a stable set of $W_0$ and therefore $\alpha(W_0)=\alpha(W)-1$.

In order to finish, it only remains to be proved that $W_0$ is a maximal induced subgraph with $\alpha(W_0)=\alpha(W)-1$.
If $u\in V(G)\setminus (V(G'_0)\cup v_1)$, then $M_u\cup M'_H$, where $M'_H$ is a maximum stable set of 
$H_{v_1}=H[V(H)\setminus (N_W(v_2)\cup v)]$ (note that $\alpha(H_{v_1})=\alpha(H)$),
is a maximum stable set of $W[V(W_0)\cup \{u\}]$  with $\alpha(W)$ vertices and
if $u\in N_W(v_2)\cap V(H)$, then $M_{v_1}\cup N_{v_1}$ ($v_2\notin M_{v_1}$ because $\{v_1,v_2\}\in E(G)$) is a maximum stable set of $W[V(W_0)\cup \{u\}]$ with $\alpha(W)$ vertices.

\medskip

\noindent $(\Rightarrow)$ 
Let us take $V_{G_0}=V(W_0)\cap V(G)$ and $V_{H_0}=V(W_0)\cap V(H)$ and let $G_0=G[V_{G_0}]$ and $H_0=H[V_{H_0}]$.
Claim 1 implies that $\alpha(G_0)\geq \alpha(G)-1$ and $\alpha(H_0)\geq \alpha(H)-1$.
Moreover, if $\alpha(G_0)= \alpha(G)-1$, then $\alpha(H_0)= \alpha(H)$ 
and if $\alpha(H_0)= \alpha(H)-1$, then $\alpha(G_0)= \alpha(G)$.
Thus, we need to considerer the following cases:

\medskip

\noindent {\bf Case $1$} If $\alpha(G_0)=\alpha(G)-1$, then there exists a maximal induced subgraph $G'_0$ with 
$\alpha(G'_0)= \alpha(G)-1$ such that $V(G_0)\subseteq V(G'_0)$ and therefore $W_0$
is as in $(i)$.

\medskip

\noindent {\bf Case $2$} If $\alpha(H_0)=\alpha(H)-1$, then there exists a maximal induced subgraph $H'_0$ with 
$\alpha(H'_0)= \alpha(H)-1$ such that $V(H_0)\subseteq V(H'_0)$ and therefore $W_0$
is as either in $(ii)$ or $(iii)$.

\medskip

\noindent {\bf Case $3$} The third case is when $\alpha(G_0)=\alpha(G)$ and $\alpha(H_0)=\alpha(H)$.
Let $M_{G_0}$ and $M_{H_0}$ be maximum stable sets of $G_0$ and $H_0$ respectively.
Using the fact that $M_{H_0}\cap \{v_1,v_2\}=\emptyset$, then we have that if $v_1,v_2\notin M_{G_0}$, 
then $M_{G_0}\cup M_{H_0}$ is a stable set of $W_0$ with $\alpha(G)+\alpha(H)$ vertices; a contradiction.
Hence $\{v_1,v_2\}\cap M_{G_0}\neq \emptyset$ for all $M_{G_0}$ maximum stable sets of $G_0$.
Furthermore, since $V_{H_0}\subset V(H)\setminus(N_W(v_k)\cup \{v\})$ if $v_k\in V_{G_0}$ 
for $k=1,2$ (there exist $M_{v_1}$ and $M_{v_2}$ stable sets with $\alpha(H)$ vertices of 
$W[V(H)\setminus(N_W(v_1)\cup \{v\})]$ and $W[V(H)\setminus(N_W(v_2)\cup \{v\})]$, respectively) 
and $\alpha(H\setminus (N_W(v_1)\cup N_W(v_2)))=\alpha(H\setminus N_H[v])=\alpha(H)-1$ 
we have that either $v_1\in M_{G_0}$ 
or $v_2\in M_{G_0}$ 
or $v_1,v_2\in M_{G_0}$
for all the maximum stable sets $M_{G_0}$ of $G_0$.%

Since $u\in M_{G_0}$ for all the maximum stable sets $M_{G_0}$ of $G_0$ if and only if 
$\alpha(G_0\setminus u)=\alpha(G_0)-1$, then either $W_0$ is as in $(i)$ whenever $v_1,v_2\in M_{G_0}$ 
for all the maximum stable sets of $G_0$%
or as in $(iv)$ in the other case. 
\qed

\begin{corollary}\label{maximalsplitt}
Let $H$ be an $\alpha$-critical graph and $v$ a vertex of $H$, then $S_0$ is a maximal 
induced subgraph of $S=s(H,v)$ with $\alpha(S_0)=\alpha(s(H,v))-1$ if and only if $S_0$ is one of the following graphs:
\[
S_0=
\begin{cases}
H\setminus v & \\
H\setminus N_{s(H,v)}(v') & \\
H\setminus N_{s(H,v)}(v'') & \\
s(H'_0,v) & \mbox{ if } v\in V(H'_0)\\
S[V(H'_0) \sqcup \{u,v'\}] & \mbox{ if } v\notin V(H'_0)\\
S[V(H'_0) \sqcup \{u,v''\}] & \mbox{ if } v\notin V(H'_0)
\end{cases}
\]
and $H'_0$ is a maximal induced subgraph of $H$ with $\alpha(H'_0)=\alpha(H)-1$.
\end{corollary}
\demo 
Since $s(H,v)=c({\cal K}_3,e,H,v)$, then the result follows directly from theorem~\ref{maximalcomposition}.
\qed

\medskip

Note that using theorem~\ref{maximalcomposition} it is easy to characterize the maximal induced subgraphs
of the edge odd-subdivision of a graph.

\medskip

Now, we will characterize the maximal induced subgraphs of the $1$-join of two graphs $G$ and $H$
with stability number equal to $\alpha(j(G,G_0,H,H_0))-1$.

\begin{theorem}\label{maximaljoin}
Let $G$ and $H$ be graphs, $G_0$ and $H_0$ be induced subgraphs of $G$ and $H$
with $\alpha(G_0)= \alpha(G)-1$ and $\alpha(H_0)= \alpha(H)-1$. 
Then $J_0$ is a maximal induced subgraph of $J=j(G,G_0,H,H_0)$ with $\alpha(J_0)=\alpha(J)-1$ if and only if
$J_0$ is one of the following graphs:
\begin{description}
\item[(i)] \hspace{1mm} $J[V(G'_0)\cup V(H'_0)]$ where $G'_0$ a maximal induced subgraph of $G$ with $\alpha(G'_0)=\alpha(G)-1$, 
$H'_0$ a maximal induced subgraph of $H$ with $\alpha(H'_0)=\alpha(H)-1$ and $\alpha(G[V(G_0)\cap V(G'_0)])=\alpha(G)-1$, 
$\alpha(H[V(H_0)\cap V(H'_0)])=\alpha(H)-1$,

\vspace{2mm}

\item[(ii)] \hspace{1mm} $J[V(G'_0) \cup V(H)]$ where $G'_0$ a maximal induced subgraph of $G$ with $\alpha(G'_0)=\alpha(G)-1$ 
and $\alpha(G[V(G_0)\cap V(G'_0)])=\alpha(G)-2$,

\vspace{2mm}

\item[(iii)] \hspace{1mm} $J[V(G) \cup V(H'_0)]$ where $H'_0$ a maximal induced subgraph of $H$ with $\alpha(H'_0)=\alpha(H)-1$ 
and $\alpha(H[V(H_0)\cap V(H'_0)])=\alpha(H)-2$,
\end{description}
\end{theorem}
\demo
$(\Leftarrow)$
$(i)$ Let $M$ be a maximum stable set of $J_0=J[V(G'_0)\cup V(H'_0)]$ and
take $M_{G'_0}=M\cap V(G'_0)$ and $M_{H'_0}=M\cap V(H'_0)$.
Since $M_{G'_0}$ and $M_{H'_0}$ are stable sets of $G'_0$ and $H'_0$ respectively, 
then $\vert M_{G'_0} \vert \leq \alpha(G)-1$ and $\vert M_{H'_0}\vert \leq \alpha(H)-1$, 
that is, $\alpha(J_0)\leq \alpha(G)+\alpha(H)-2=\alpha(J)-1$.
Moreover, if $M'_G$ and $M'_H$ are maximum stable sets of $J[V(G'_0)\cap V(G_0)]$ and $J[V(H'_0)\cap V(H_0)]$ respectively, 
then $M=M'_G\cup M'_H$ is a stable set of $J_0$ with $\alpha(J)-1$ vertices
and therefore $\alpha(J_0)=\alpha(J)-1$.

Now, let $v\in V(J)\setminus (V(G'_0)\cup V(H'_0))$ and assume that $v\in V(G)\setminus V(G'_0)$. 
By the maximality of $G'_0$ we have that $\alpha(J[V(G'_0)\cup \{v\}])=\alpha(G)$.
Since $J[V(G'_0)\cup \{v\}]$ and $H[V(H_0)\cap V(H'_0)]$ are not connected by some edge of $J$, 
then $\alpha(J[V(G'_0)\cup \{v\}\cup v(H'_0)])=\alpha(G)+\alpha(H[V(H_0)\cap V(H'_0)])=\alpha(J)$.
Using the same arguments for $v\in V(H)\setminus V(H_0)$ we can conclude 
that $J_0$ is a maximal induced subgraph of $J$ with $\alpha(J_0)=\alpha(J)-1$.

\medskip

$(ii)$ Let $M$ be a maximum stable set of $J_0=J[V(G'_0)\cup V(H)]$ and take $M_{G'_0}=M\cap V(G'_0)$ and $M_H=M\cap V(H)$.
Since $M_{G'_0}$ and $M_H$ are stable sets of $G'_0$ and $H$ respectively, 
then $\vert M_{G'_0} \vert \leq \alpha(G)-1$ and $\vert M_H\vert \leq \alpha(H)$. 
Furthermore, if $\vert M_H \vert=\alpha(H)$, then $M_H\cap V(H)\setminus V(H_0)\neq \emptyset$ and $M_{G'_0}\subset V(G_0)\cup V(G'_0)$. 
Thus, $\vert M_{G'_0}\vert \leq \alpha(G)-2$ and therefore $\alpha(J_0)\leq \alpha(J)-1$.

Now, since $\alpha(G[V(G_0)\cap V(G'_0)])=\alpha(G)-2$, then there exists $M'$ a stable set of $G'_0$ with $\alpha(G)-2$ vertices.
Taking $M=M'\cup M_H$, where $M_H$ is a maximum stable set of $H$, we can conclude that 
$M$ is a stable set of $J_0$ with $\alpha(G)+\alpha(H)-2$ vertices and therefore $\alpha(J_0)= \alpha(J)-1$.

Moreover, $J_0$ is maximal with $\alpha(J_0)=\alpha(J)-1$ because if $v\in V(G)\setminus V(G'_0)$, 
then by the maximality of $G'_0$ we have that there exists a stable set $M_v$ of $G$ such that
$v\in M_v\subset V(G'_0)\cup \{v\}$ and $\vert M_v\vert=\alpha(G)$.
Thus, $M=M_v\cup M_{H}$ (where $M_H$ is a maximum stable set of $H$) is a stable set of 
$J[V(G'_0)\cup \{v\}\cup V(H)]$ with $\vert M\vert=\vert M_v\vert+\vert M_H\vert=\alpha(G)+\alpha(H)-1=\alpha(J)$
vertices.

$(iii)$ The result is followed using exactly the same argument as in $(ii)$.

\medskip

\noindent $(\Rightarrow)$
Let $V_G=V(J_0)\cap V(G)$, $V_H=V(J_0)\cap V(H)$, $G'=G[V_G]$ and $H'=H[V_H]$.
Since $G'_0$ and $H'_0$ are induced subgraphs of $G$ and $H$ respectively,
then $\alpha(G'_0)\leq \alpha(G)$ and $\alpha(H'_0)\leq \alpha(H)$. 
Furthermore, is not difficult to see that $\alpha(G)+\alpha(H)-2\leq \alpha(G')+\alpha(H')\leq \alpha(G)+\alpha(H)-1$.
Thus, $\alpha(G_0)\geq \alpha(G)-2$ and $\alpha(H'_0)\geq \alpha(H)-2$; moreover if $\alpha(G'_0)=\alpha(G)-k$ for some $k=0,1,2$, then $\alpha(H'_0)\geq \alpha(H)-k$. 

Now, we will considerer the different options for the values of $\alpha(G'_0)$.
If $\alpha(G'_0)=\alpha(G)$, then $\alpha(H'_0)\geq \alpha(H)-2$ and therefore $J_0$ is as in $(iii)$.
If $\alpha(G'_0)\leq \alpha(G)-1$, then $\alpha(H'_0)\geq \alpha(H)-1$,
that is, $\alpha(H)-1\leq \alpha(H'_0)\leq\alpha(H)$.
Finally, if $\alpha(H'_0)=\alpha(H)-1$, then $J_0$ is as in $(i)$ and
if $\alpha(H'_0)=\alpha(H)$, then $J_0$ is as in $(ii)$.
\qed

\begin{corollary}
Let $G$ be a graph and $v$ a vertex of $G$, then $D_0$ is a maximal induced 
subgraph of $d(G,v)$ with $\alpha(D_0)=\alpha(d(G,v))-1$ if and only if
$D_0$ is one of the following graphs:
\[
D_0=
\begin{cases}
G'_0 & \mbox{ if } G'_0 \mbox{ is a maximal induced subgraph of } G \mbox{ with } \alpha(G'_0)=\alpha(G)-1\\
& \mbox{ and } \alpha(G[V(G_0)\setminus N[v]])=\alpha(G)-1\\
G[V(G'_0)\cup \{v\}] & \mbox{ if } G´'_0 \mbox{ is a maximal induced subgraph of } G \mbox{ with } \alpha(G'_0)=\alpha(G)-1\\
& \mbox{ and } \alpha(G[V(G'_0)\setminus N[v]])=\alpha(G)-2.
\end{cases}
\]
\end{corollary}
\demo 
Since $d(G,v)=j(G,G\setminus N[v],{\cal K}_1,\emptyset)$, then the result is followed directly from theorem~\ref{maximaljoin} .
\qed


\section{Basic $\alpha$-critical graphs}

An $\alpha$-critical graph $G$ is called {\it basic} if it is splitting free (note that if $G$ is a splitting free graph, then $G$ is odd subdivision free) and duplication free.
A basic $\alpha$-critical graph $G$ is called {\it strongly basic} if it is edge-vertex composition free and $1$-join composition free. 
Lov\'asz and Plummer defined in \cite[pag. 453]{LovaszPlummer} that a $\tau$-critical graph
is basic if it is splitting free. 

\medskip

By \cite[Lemmma 12.1.4]{LovaszPlummer} we have that a connected $\alpha$-critical graph $G$ is splitting free if and only if $deg(v)\geq 3$ for all $v\in V(G)$.
Therefore is clear that if $G$ and $H$ are splitting free, then $j(G,G_0,H,H_0)$ is splitting free. 
However is possible that $G$ or $H$ are not splitting free and $j(G,G_0,H,H_0)$ would be splitting free

The next theorem characterize when the $1$-join of two graphs is a basic $\alpha$-critical graph.
 
\begin{theorem}
Let $G$ and $H$ be graphs, $G_0$ and $H_0$ be induced graphs of $G$ and $H$
with $\alpha(G_0)= \alpha(G)-1$ and $\alpha(H_0)= \alpha(H)-1$ and let $J=j(G,G_0,H,H_0)$.
Then
\begin{description}
\item[{\it (i)}] \hspace{1mm} $J$ is a connected $\alpha$-critical splitting free graph if and only if
the following conditions are satisfied:
\begin{description}
\item[$\bullet$] \hspace{1mm} $G_0$ and $H_0$ are connected, 

\vspace{2mm}

\item[$\bullet$] \hspace{1mm} if $v\in V(G_0)$, then ${\rm deg}_G(v)\geq 3$ and if $v\in V(H_0)$, then ${\rm deg}_H(v)\geq 3$,
in particular, we have that $G_0$ and $H_0$ are splitting free graphs.

\vspace{2mm}

\item[$\bullet$] \hspace{1mm} if either $V(G_0),V(H_0)\neq \emptyset$ or $V(G_0)=\emptyset$, 
then either $H\neq s(H',v)\setminus u$ for any graph $H'$ and for all $v\in V(H')$
or $H_0\neq H'\setminus v$ for any $v\in V(H')$ such that $H=s(H',v)\setminus u$ for some $H'$,

\vspace{2mm}

\item[$\bullet$] \hspace{1mm} $G$, $G_0$ and $H$, $H_0$ satisfy the properties $(i),(ii)$ and $(iii)$ in theorem~\ref{1join}.
\end{description}

\hspace{1mm}

\item[{\it (ii)}] \hspace{1mm} $J$ is a connected $\alpha$-critical odd subdivision free graph if and only if 
$G$ and $H$ are connected, $G_0$ and $H_0$ are odd subdivision free graphs and $G$, $G_0$ and $H$, $H_0$ 
satisfy the properties $(i),(ii)$ and $(iii)$ in theorem~\ref{1join}.

\hspace{1mm}
\item[{\it (iii)}] \hspace{1mm} $J$ is a duplication free graph if and only if 
$G$ and $H$ are duplication free graphs and either $G_1=j(G,G_0,{\cal K}_1,\emptyset)$ is duplication free ($G_0\neq G\setminus N[u]$ for all $u\in V(G)$) 
or $H_1=j(H,H_0,{\cal K}_1,\emptyset)$ is duplication free ($H_0\neq H\setminus N[v]$ for all $v\in V(H)$).
\end{description}
\end{theorem}
\demo
$(i)$ ($\Rightarrow$)
If $v\in V(G_0)$, then $s(j(G,G_0,H,H_0),v)=j(s(G,v),s(G_0,v),H,H_0)$.
Hence, if we assume that $G_0$ or $H_0$ are not splitting-free, 
then $j(G,G_0,H,H_0)$ is not splitting-free; a contradiction.
Moreover, if $G=s(G',v)$ for some graph $G'$ and some $v\in V(G')$ and $G_0$ contain the vertex $u$ of degree two 
added in the splitting of $v$ in $G'$, then $j(G,G_0,H,H_0)$ is not splitting-free; a contradiction.

Is not difficult to see 
that if either $G_0$ or $H_0$ are not connected, then $j(G,G_0,H,H_0)$ will not be connected; a contradiction.
Finally, if $V(G_0)=\emptyset$, $H=s(H',v)\setminus u$ for some graph $H'$, $v\in V(H')$
and $H'_0=H'\setminus v$,
then $j(G,G_0,H,H_0)=s(H',v)$; a contradiction.

\medskip

\noindent ($\Leftarrow$)
Let us assume that $j(G,G_0,H,H_0)$ is not splitting-free, then $j(G,G_0,H,H_0)$ has a vertex $v$ of degree two. 
Moreover, we can assume without lost of generalization that $v$ is in $V(G)$.
Let $v$ be a vertex of $G$, then ${\rm deg}_{j(G,G_0,H,H_0)}(v)=2$ if and only if
\[
{\rm deg}_{G}(v)=
\begin{cases}
2 & \mbox{ if } v\in V(G_0),\\
1 & \mbox{ if } v\in V(G)\setminus V(G_0) \mbox{ and } \vert V(H)\setminus V(H_0)\vert=1, \\
0 &  \mbox{ if } v\in V(G)\setminus V(G_0) \mbox{ and } \vert V(H)\setminus V(H_0)\vert=2.
\end{cases}
\]
Then we need to considerer the following three cases:
\begin{description}
\item[$\bullet$] ${\rm deg}(v)=2$ and $v\in V(G_0)$.
Then, since $N(v)\subset V(G)$ we will have that $G$ is not splitting-free
and therefore $G_0$ contains a vertex of degree two; 
a contradiction, 

\vspace{2mm}

\item[$\bullet$] ${\rm deg}(v)=1$, $\vert V(H)\setminus V(H_0)\vert=1$ and $v\in V(G)\setminus V(G_0)$.
Let $u'=V(H)\setminus V(H_0)$ and $e=\{u',w\}$ any edge in $H$ incident with $u'$.
Since $\alpha(H_0)=\alpha(H)-1$, then $e$ is not $\alpha$-critical edge of $H$; a contradiction. 

\vspace{2mm}

\item[$\bullet$] ${\rm deg}(v)=0$, $\vert V(H)\setminus V(H_0)\vert=2$ and $v\in V(G)\setminus V(G_0)$.
In this case we have that if $V(G_0)\neq \emptyset$, then $j(G,G_0,H,H_0)$ would not be connected
and if $V(G_0)= \emptyset$, then $H=s(G',w)\setminus w$ for some graph $G'$ and some $w\in V(G')$; a contradiction.
\end{description}

Finally, it is clear that if $G_0$ and $H_0$ are connected, then $j(G,G_0,H,H_0)$ is also connected.

\medskip

$(ii)$ ($\Rightarrow$) If $e=\{x,y\}\in E(G)$, then 
\[
s(j(G,G_0,H,H_0),e)=
\begin{cases}
j(s(G,e),s(G_0,e),H,H_0) & \mbox{ if } e\in E(G_0),\\
j(s(G,e),G[V(G_0)\sqcup \{x',y'\}],H,H_0) & \mbox{ if } e\in E(G)\setminus E(G_0).
\end{cases}
\]
Note that there exists a similar relation when $e\in E(H)$.
Hence we can conclude that if either $G_0$ or $H_0$ are not odd subdivision-free, then $j(G,G_0,H,H_0)$ is not odd subdivision-free.

\medskip

For the converse, assume that $j(G,G_0,H,H_0)$ is not odd subdivision-free, 
that is, there exist $x'$ and $y'$ two adjacent vertices with degree two with non-adjacent neighborhoods.
Hence using the observation made at the beginning of the proof we need to consider the following possibilities:
\begin{description}
\item[$\bullet$] $x',y'\in V(G_0)\cup V(H_0)$. 
Since $x'$ and $y'$ are adjacent we have that either $x',y'\in V(G_0)$ or $x',y'\in V(H_0)$
and therefore either $G$ or $H$ are not odd subdivision-free.

\vspace{2mm}

\item[$\bullet$] $x'\in V(G_0)$ and $y'\in V(G)\setminus V(G_0)$. 
Let $y=N(y')\setminus x'=V(H)\setminus V(H_0)$.
Since $\alpha(H_0)=\alpha(H)-1$, then all the edges $e$ of $H$
incidents with $y$ are not $\alpha$-critical edges.
Therefore by theorem~\ref{1join} we have that $j(G,G_0,H,H_0)$ can not be $\alpha$-critical; a contradiction.

\vspace{2mm}

\item[$\bullet$] $x',y'=V(G)\setminus V(G_0)$.
In this case we have that $w=V(H)\setminus V(H_0)$ is adjacent to $x'$ and $y'$ and therefore $\{x',y',w\}$
induce a triangle (a connected component); a contradiction to the connectivity of $j(G,G_0,H,H_0)$.

\vspace{2mm}

\item[$\bullet$] $x'\in V(H_0)$ and $y'\in V(H)\setminus V(H_0)$ or $x',y'=V(H)\setminus V(H_0)$. 
Exactly the same situation that in the previous cases.

\vspace{2mm}

\item[$\bullet$] $x'=V(G)\setminus V(G_0)$ and $y'=V(H)\setminus V(H_0)$.
Let $x=N(x')\setminus y'\in V(G_0)$ and $y=N(y')\setminus x'\in V(H_0)$.
Since $\alpha(H_0)=\alpha(H)-1$, then the edge $e_1=\{x,x'\}$ is not an $\alpha$-critical edge of $G$
and by theorem~\ref{1join} we have that it is a contradiction to the fact that $j(G,G_0,H,H_0)$ is an $\alpha$-critical graph.
\end{description}

$(iii)$ ($\Rightarrow$)
Let $v \in V(G)$, then
\[
d(j(G,G_0,H,H_0),u)=
\begin{cases}
j(d(G,u),d(G_0,u),H,H_0) & \mbox{ if } u\in V(G_0),\\
j(d(G,u),G_0,H,H_0) & \mbox{ if } u\in V(G)\setminus V(G_0).
\end{cases}
\]
Also there exists a similar relation when $u\in V(H)$. 
Hence we can conclude that if either $G$ or $H$ is not duplication-free, then $j(G,G_0,H,H_0)$ is not duplication-free.

To finish this part, is not difficult to realized that if $G_0=G\setminus N[u]$ for some $u\in V(G)$ 
and $H_0=H\setminus N[v]$ for some $v\in V(H)$, then 
$$
j(G,G_0,H,H_0)=d(j(G\setminus u,G_0\setminus u,H,H_0),v)=d(j(G,G_0,H\setminus v,H_0\setminus v),u),
$$
that is, $j(G,G_0,H,H_0)$ is not duplication-free; a contradiction.

\medskip

\noindent ($\Leftarrow$)
Assume that $j(G,G_0,H,H_0)$ is not duplication-free, that is, that there exist $u,v$ vertices of $j(G,G_0,H,H_0)$
such that $u$ is the duplication of $v$.

We can assume without lost of generalization that $u\in V(G)$.
If $u\in V(G_0)$, then $v\in V(G_0)$ and therefore $u$ is the duplication
of $v$ in $G$, that is, $G$ is not duplication-free; a contradiction.

If $u \in V(G)\setminus V(G_0)$, then we have that either $v \in V(G)\setminus V(G_0)$ or $v \in V(H)\setminus V(H_0)$.
In the first case we have that $u$ is the duplication of $v$ in $G$ and 
in the second case we have that $N_G[u]=V(G)\setminus V(G_0)$ and $N_H[v]=V(H)\setminus V(H_0)$ 
because $N[u]=N[v]=(V(G)\setminus V(G_0))\cup (V(H)\setminus V(H_0))$ in $j(G,G_0,H,H_0)$.
Therefore $G_0=G\setminus N[u]$ and $H_0=H\setminus N[v]$; a contradiction.
\qed



\end{document}